\documentclass[12pt]{article}
\usepackage{amssymb}
\usepackage{epsfig}
\makeatletter

\@addtoreset{equation}{section}

\makeatother

\newtheorem{thm}[equation]{Theorem}
\newtheorem{lemma}[equation]{Lemma}

\renewcommand\thefigure{\thesection.\@arabic\c@figure}

\renewcommand\thetable{\thesection.\@arabic\c@table}

\def\reff#1{(\ref{#1})}
\def\sobre#1#2{\lower 1ex \hbox{ $#1 \atop #2 $ } }

\begin{document}

\def\E{{\mathbb E}}
\def\P{{\mathbb P}}
\def\R{{\mathbb R}}
\def\Z{{\mathbb Z}}
\def\V{{\mathbb V}}
\def\N{{\mathbb N}}
\def\NN{{\bf N}}
\def\X{{\cal X}}
\def\Y{{\bf Y}}
\def\G{{\cal G}}
\def\T{{\cal T}}
\def\C{{\bf C}}
\def\D{{\bf D}}
\def\G{{\bf G}}
\def\U{{\bf U}}
\def\K{{\bf K}}
\def\H{{\bf H}}
\def\n{{\bf n}}
\def\b{{\bf b}}
\def\g{{\bf g}}
\def\mm{{m}}
\def\sqr{\vcenter{
         \hrule height.1mm
         \hbox{\vrule width.1mm height2.2mm\kern2.18mm\vrule width.1mm}
         \hrule height.1mm}}                  
\def\square{\ifmmode\sqr\else{$\sqr$}\fi}
\def\one{{\bf 1}\hskip-.5mm}
\def\liml{\lim_{L\to\infty}}
\def\given{\ \vert \ }
\def\a{\alpha}
\def\ze{{\zeta}}
\def\be{{\beta}}
\def\de{{\delta}}
\def\e{{\epsilon}}
\def\la{{\lambda}}
\def\ga{{\gamma}}
\def\th{{\theta}}
\def\proof{\noindent{\bf Proof. }}
\def\sign{\hbox{sign}}
\def\mod{\hbox{mod}}
\def\rate{{e^{- \beta|\ga|}}}
\def\A{{\bf A}}
\def\B{{\bf B}}
\def\C{{\bf C}}
\def\D{{\bf D}}

\def\={&=&}
\def\+{&+&}

 \title{A dynamic one-dimensional interface\\ 
interacting with a wall}
 \date{}
      
 \author{F. M. Dunlop, P. A. Ferrari, L. R. G. Fontes}

 \maketitle

 \noindent{\bf Abstract.} 
 We study a symmetric randomly moving line interacting by exclusion with a
 wall. We show that the expectation of the position of the line at the origin
 when it starts attached to the wall satisfies the following bounds:
 $$
 c_1t^{1/4} \le\E\xi_t(0) \le c_2 t^{1/4}\log t
 $$
 The result is obtained by comparison with a ``free'' process, a random
 line that has the same behavior but does not see the wall.  The free process
 is isomorphic to the symmetric nearest neigbor one-dimensional simple
 exclusion process. The height at the origin in the interface model
 corresponds to the integrated flux of particles through the origin in the
 simple exclusion process.  We compute explicitly the asymptotic variance of
 the flux and show that the probability that this flux exceeds $Kt^{1/4}\log
 t$ is bounded above by const.~$t^{2-K}$.  We have also performed numerical
 simulations, which indicate $\E\xi_t(0)^2 \sim t^{1/2}\log t$ as
 $t\to\infty$.  \bigskip
 
 

\section{Introduction}

We consider a process $\xi_t$ on 
$$
\X=\{\xi\in\N^\Z: |\xi(x)- \xi(x+1)|=1,\, \xi(0) \hbox{ even}\} 
$$
the space of trajectories of nearest neighbor random walks that stay non
negative and such that at even ``times'' the walk visits even integers.

The generator of the process is given by
\begin{equation}
  \label{gen}
  {\cal L}f(\xi) = \frac12 \sum_x \one\{\xi + \Delta
  \xi(x) \,\delta_x \ge 0\}\, [f(\xi + \Delta
  \xi(x)\,\delta_x) - f(\xi)]
\end{equation}
where $\delta_x$ is the infinite vector having $1$ in the $x^{\rm th}$
coordinate and zero on the others. The sum $\xi+ a\delta_x$ is
understood coordinatewise. The discrete Laplacian $\Delta$ is defined
by 
\begin{equation}
  \label{p28}
  \Delta \xi(x) := \xi(x+1) - 2 \xi(x) + \xi(x-1).
\end{equation}
In words we can describe the dynamics as follows. The discrete Laplacian
assumes only three values, $-2$, $0$ and $2$. When the Laplacian is zero, the
interface does not move. When it is $-2$ or $2$, at rate $\frac12$ it makes a
jump of length $2$ in the same direction as the sign of the Laplacian. Over
this motion we impose a restriction to keep the process in $\X$: the interface
cannot be negative, so we simply prohibit the jumps which violate the
restriction. This is the meaning of the indicator function $\one\{\xi + \Delta
\xi(x) \,\delta_x \ge 0\}$ in the generator.  We can think the prohibition of
becoming negative as the interaction by exclusion of the interface with a wall
at $-1$. For shortness we call $\xi_t$ the \emph{wall process}.

Our main result is the following

\begin{thm}\label{98}
  Let $\xi_t$ be the process with generator \reff{gen} and initial \emph{flat
    configuration}:
\begin{equation}
  \label{99}
   \xi_0(x) := x\,(\hbox {\rm mod } 2). 
\end{equation}
Then there exist positive constants $c_1$ and $c_2$ such that 
\begin{equation}
  \label{100}
   c_1t^{1/4} \le\E\xi_t(0) \le c_2 t^{1/4}\log t
\end{equation}
for sufficiently large $t$.
\end{thm}

Theorem \ref{98} catches the effect of the ``entropic repulsion'' in a
stochastically moving interface interacting with a wall by exclusion. 

The line induced in $\R^2$ by joining 
$(x,\xi_t(x))$ to $(x+1,\xi_t(x+1))$ for all $x\in \Z$
 has the same behavior as the
interface between $-1$'s and $1$'s in a zero-temperature two-dimensional
nearest-neighbors Ising model with a positive external field in the semiplane
below the diagonal $x=y$ and with initial condition ``all plus'' below the
diagonal and ``all minus'' above it. See Section \ref{ising} for details.

The equilibrium statistical mechanics of this model is well known. If one
considers the generator ${\cal L}$ restricted to the box $[-L,L]$ with
boundary conditions $\xi_t(-L)=\xi_t(L)=0$, the invariant distribution is the
uniform distribution in the set of nearest neighbors random walk trajectories
starting at time $-L$ at the origin, finishing at time $L$ at the origin and
being non negative for all intermediate times. Actually the uniform measure is
even reversible for the process. But the uniform measure in this set
corresponds to the law of a symmetric nearest neighbors random walk $X_i$
conditioned to the set $\{X_{-L}=X_L=0,\, X_i\ge 0, i\in [-L,L]\}$.  Hence,
the typical height of a configuration $\xi$ with the invariant law in the bulk
of the box is
\[
\xi([rL]) \;\sim\; O(\sqrt{rL})
\]
More precisely, the normalized process process $(L^{1/2}\xi([Lr]),
r\in[-1,1])$ converges as $L\to\infty$ to Brownian excursion on $[-1,1]$; see
Theorem 2.6 of Kaigh (1976).

Many papers deal with the problem of entropic repulsion in Equilibrium
Statistical Mechanics.  The role of the entropic repulsion in the Gaussian
free field was studied by 
Lebowitz and Maes (1987), 
Bolthausen, 
Deuschel and Zeitouni (1995), 
Deuschel (1996) and 
Deuschel and Giacomin (1999).

Entropic repulsion for Ising, SOS and related models was discussed in
Bricmont, El Mellouki and Fr\"ohlich (1986),
Bricmont (1990),
Holick\'y and Zahradn\'\i k (1993), 
Cesi and Martinelli (1996), 
Lebowitz, Mazel and Suhov (1996),
Dinaburg and Mazel (1994) and 
Ferrari and Mart\'{\i}nez (1998). 

The exponent $1/4$ for dynamic entropic repulsion was predicted by Lipowsky
(1985) using scaling arguments. This exponent was then found numerically
by 
Mon, Binder, Landau (1987),
Albano, Binder, Heermann, Paul (1989-1992), see Binder (1990),
De Coninck, Dunlop and Menu (1993).
It has also been observed in real experiments by
Bartelt, Goldberg, Einstein, Williams, Heyraud, M\'etois (1993).
Further theoretical investigations of dynamics of lines, in relation to 
experiments can be found in Blagojevic, Duxbury (1999).

Dynamic entropic repulsion for a line of finite extension $L$ when 
$t,L\to\infty$ strongly depends on the ratio $t/L^2$. The present paper deals 
with $L=\infty$ (analytical) or $t/L^2\to 0$ (numerical). The case 
$t/L^2={\cal O}(1)$ has been studied by Funaki and Olla (2001).

The exponent $1/4$ also applies to the growth of fluctuations of an initially
straight interface not interacting with the wall (see \reff{105p} below). For
the Gaussian case, explicit computations were made by Abraham, Upton (1989),
Abraham, Collet, De Coninck, Dunlop (1990). It was observed numerically in the
two-dimensional Ising model by Stauffer, Landau (1989).

The strategy to show Theorem \ref{98} is to compare the wall process $\xi_t$
with a \emph {free} process $\zeta_t$ having the same local
dynamics as $\xi_t$ but not interacting with the wall. The free process
lives in
$$
\X_0 = \{\zeta\in \Z^\Z: |\zeta(x)-\zeta(x+1)|=1, \zeta(0)= \hbox{ even}\}
$$
and its generator is 
\begin{equation}
  \label{gen0}
  {\cal L}_0f(\zeta)\; =\; \frac12\sum_x  [f(\zeta + \Delta
  \zeta(x)\,\delta_x) - f(\zeta)]
\end{equation}

In the next section we prove that with flat initial condition the variance
of the height at the origin for the free process behaves as $t^{1/2}$:
\begin{equation}
 \label{vfp}
 \lim_{t\to\infty}t^{-1/2}\V \zeta_t(0)
 \;=\;\frac{1}{\sqrt\pi}.
\end{equation}

We then couple the wall process and the free process in such a way that
\begin{equation}
  \label{cou}
  \zeta_t(x) \le \xi_t(x)
\end{equation}
for all $x$ and $t$. 
The free process has enough symmetry and, properly rescaled, has uniformly
bounded in time exponential moments.
With these, \reff{vfp} and \reff{cou}, we get the lower bound in
\reff{100}.

The idea for the upperbound is to consider a family of free processes
with initial condition depending on $t$:
$$
\zeta^{a_t}_0 (x) = \zeta_0(x) + a_t 
$$
(a flat interface of height $a_t$). Then we fix $a_t = c t^{1/4}\log t$,
the constant $c$ to be determined later and exhibit a coupling under which 
\begin{equation}
  \label{eq:c01}
  \xi_s(0) \le \zeta^{a_t}_s(0)
\end{equation}
for all $s\le t$ with large probability. Combined with \reff{vfp}, inequality
\reff{eq:c01} is the key for the upperbound in \reff{100}. The existence of
exponential moments (mentioned above) yields the moderate deviations result
needed here.

The control of the fluctuations of the position at the origin of the free
process is obtained by an isomorphism between the free process and the
one-dimensional symmetric nearest-neighbor simple exclusion process. Under
this map, $\zeta_t(0)=2J_t$, where $J_t$ is the integrated flux of particles at
the origin for the exclusion process. We compute explicitly the asymptotic
variance of the integrated flux for the flat
initial condition in Theorem \ref{104} below and obtain
\begin{equation}
  \label{vfp6}
  \lim_{t\to\infty}{\V J_t \over \sqrt t} 
 =\frac{1}{4\sqrt{\pi}}.
\end{equation}
De Masi and Ferrari (1985) proved that the asymptotic variance of the
integrated flux when the initial configuration is distributed according to a
product measure with density $1/2$ is given by
\begin{equation}
  \label{105a}
  \lim_{t\to\infty}{\V J_t \over \sqrt t} 
 =\frac{1}{2\sqrt{2\pi}}.
\end{equation}
which is strictly bigger than \reff{vfp6}. When the initial density is $\rho$,
the asymptotic variance is given by $\rho(1-\rho)\sqrt{2/\pi}$.  The method to
show \reff{vfp6} and \reff{105a} is based on duality and comparison with
systems of independent particles and it is inspired by Arratia (1983), who
used these tools to compute the variance of a tagged particle for the process
starting with an (invariant) product measure. However a modification of
Arratia's proof is needed in \reff{vfp6} due to the deterministic character of
the initial configuration.

The study of the flux in the simple exclusion process is done in Section 2. In
Section 3 we prove Theorem \ref{98}. Section 4 is devoted to numerical 
simulations.
Section 5 shows the equivalence to the dynamics of a particular zero 
temperature Ising model interface.

\section{ The free process and simple exclusion}

The \emph{simple exclusion process} lives in $\{0,1\}^\Z$ and its generator is
\begin{equation}
  \label{101}
  {\cal L}^{ex}f(\eta) = \frac12\sum_{x\in\Z} [f(\eta^{x,x+1}) - f(\eta)]
\end{equation}
where
\begin{equation}
  \label{130}
  \eta^{x,x+1}(y) = \cases{\eta(y) &if $y\neq x, x+1$\cr
\eta(x+1) &if $y= x$\cr
\eta(x) &if $y= x+1$\cr}
\end{equation}
It is convenient to construct the processes using the Harris graphical
construction.

\medskip
\noindent{\bf Harris graphical construction}
Let $(N_t(x):x\in\Z)$ be a family of independent Poisson processes of rate
$\frac12$. For each $x$, $N_t(x)$ counts the number of Poisson events
associated to $x$ in the time interval $[0,t]$.  Denote $d N_t(x) =
\one\{$there is a Poisson event associated to $x$ at time $t\}=\lim_{h\to0}
(N_{t}(x) - N_{t-h}(x))$. Let $\eta_t$ be the process defined by
\begin{equation}
  \label{111}
  d \eta_t(x) = (\eta_t(x-1) - \eta_t(x)) dN_t(x) + (\eta_t(x+1) -
  \eta_t(x)) dN_t(x+1).  
\end{equation}
The process is well defined because for each finite time $t$ the value of the
process in a finite box can be determined by looking at only a finite but
random number of Poisson events and initial values. See for instance Arratia
(1983). In words, the motion can be described as follows. The Poisson marks of
$N_t(x)$ are associated to the bond $(x,x+1)$ and each time a Poisson mark
occurs, the contents of the associated bond are interchanged. It is immediate
to show that this process has generator \reff{101}.

\medskip
\noindent{\bf Stirring particles.} 
To introduce the notion of duality and to deal with the flux of particles it
is convenient to follow the ``stirring particles'' as
defined by Arratia (1983). Let $X^x_t$ be
the position at time $t$ determined by $X^x_0=x$ and the equations
\begin{equation}
  \label{117}
  d X^x_t =  dN_t(X^x_t)- dN_t(X^x_t-1)
\end{equation}
So that, each time a Poisson mark associated to one of the neighboring bonds
of a particle occurs, the particle jumps across the bond. Of course, if both
extremes of a bond are occupied, the particles jump simultaneously, respecting
the exclusion condition ``at most one particle per site''. For each $t\ge 0$
the (random) map
\begin{equation}
  \label{119}
  x\mapsto X^x_t
\end{equation}
is a bijection of $\Z$ in $\Z$. The (marginal) law of $X^x_t$ is a symmetric
nearest neighbor random walk starting at $x$. 

\medskip
\noindent{\bf Duality.}
Let $y\mapsto D^y_t$ be the inverse map defined by $x = D^y_t$ if and only if $y =
X^x_t$.  The following ``duality formula'' holds immediately
\begin{equation}
  \label{118}
  \eta_t(y) = \eta_0(D^y_t) 
\end{equation}
So, 
\begin{equation}
  \label{120}
  \prod_{y\in A}\eta_t(y) = \prod_{y\in A}\eta_0(D^y_t) 
\end{equation}
Notice that for a finite set of sites $A$, $\{D^y_t:y\in A\}$ has the same
one-time marginal as a simple exclusion process with initial condition $A$
(here we are identifying the configuration $\eta$ with the set
$\{x:\eta(x)=1\}$). When $A=\{y\}$ (contains only one site), the one-time
marginal $D^y_t$ has the same law as $X^y_t$ for all $t\ge 0$.

\paragraph{Integration by parts formula}
Consider $(Y^i_t,Z^j_t)$ independent random walks with the same marginals as
the stirring process $(X^i_t,X^j_t)$. The generator of the process
$(Y^i_t,Z^j_t)$ is the following:
\begin{equation}
  \label{110a}
  Uf(i,j) = \frac12 \sum_{e\in\Z^2:|e|=1}[f((i,j)+e) -  f(i,j)] 
\end{equation}
and the generator of the process $X_t$ is
\begin{equation}
  \label{110b}
  Vf(i,j) = \cases{\frac12 \sum_{e\in\Z^2:|e|=1}[f((i,j)+e) - f(i,j)] &if
    $i-j>1$\cr 
\frac12 f(j,i) + \frac12 f(i+1,j)+ \frac12 f(i,j-1)- \frac32 f(i,j)
&if 
$i-j=-1$ \cr
\frac12 f(j,i) + \frac12 f(i-1,j)+ \frac12 f(i,j+1)- \frac32 f(i,j)
&if 
$i-j=1$ }
\end{equation}
Hence, for $i\neq j$,  
\begin{equation}
  \label{110c}
  Uf(i,j)- Vf(i,j) = -\,\frac12\one\{|i-j|=1\} (f(i,j)
  +f(j,i) 
  -f(i,i) -f(j,j))  
\end{equation}
Let $U_t$ and $V_t$ be the semigroups generated by $U$ and $V$ respectively.
Let $f:\Z^2\to\R$. Then
\begin{eqnarray}
  \label{110d}
  \E f(X^i_t,X^j_t) - \E f(Y^i_t,Z^j_t) = [V_t - U_t] f(i,j)
  \;=\; \int_0^t V_s [U-V] U_{t-s}f(i,j) ds
\end{eqnarray}
where the last identity is the integration by parts formula (see Liggett
(1985) Proposition 8.1.7). Now, using \reff{110c} to compute \reff{110d} we
get for $i\neq j$:
\begin{eqnarray}
  \label{110e}
\lefteqn{\E f(X^i_t,X^j_t) - \E f(Y^i_t,Z^j_t) }\nonumber\\
&=&-\frac12 \int_0^t \,ds\;
\E\Bigl(\one\{|X^i_s-X^j_s|=1\}\\
&&\qquad\times\,
[f(Y^{X^{i}_{t-s}}_t,Z^{X^{j}_{t-s}}_t)+
  f(Y^{X^{j}_{t-s}}_t,Z^{X^{i}_{t-s}}_t) -
  f(Y^{X^{i}_{t-s}}_t,Z^{X^{i}_{t-s}}_t) 
  - f(Y^{X^{j}_{t-s}}_t,Z^{X^{j}_{t-s}}_t)]\Bigr).\nonumber
\end{eqnarray}
This identity will be used in the sequel.
\medskip

\noindent{\bf Flux.}
Let $J_t$ be the integrated flux of $\eta$ particles through the point $-1/2$ 
in the exclusion process:
\begin{equation}
  \label{eq:1}
  J_t := \sum_{x< 0} \eta_0(x)\one\{X^{x}_t\ge 0\}\, - \,\sum_{x\ge 0}
  \eta_0(x)\one\{X^{x}_t< 0\}  
\end{equation}
where $X^{x}_t$ is the position at time $t$ of the exclusion particle that at
time zero was at position $x$. 

\medskip

Replacing \reff{99} in \reff{eq:1}, we write
\begin{equation}
  \label{115}
  J_t := \sum_{i< 0} \one\{X^{2i}_t\ge 0\}\, - \,\sum_{i\ge 0} \one\{X^{2i}_t<
  0\}.  
\end{equation}

$J_t$ is almost symmetric. Let 
\begin{equation}
  \label{1101}
  H_t := \sum_{i<0} \one\{X^{2i}_t\ge 0\};\,\,\,
  H'_t := \,\sum_{i\ge 0} \one\{X^{2i}_t<-1\};\,\,\, 
  I_t := \,\sum_{i\ge 0} \one\{X^{2i}_t<0\}.  
\end{equation}
Then, clearly, 
\begin{equation}
  \label{1102}
  J_t= H_t-I_t;\,\,\,H_t\sim H'_t;\,\,\,|H'_t-I_t|\leq1,
\end{equation}
where $\sim$ means identity in distribution and is justified in this case by
spatial and distributional symmetry.
\begin{thm}
  \label{104}
For the simple exclusion process with generator \reff{101} and initial
condition $\eta_0 = 1-2\xi_0$, as defined in \reff{99},  
\begin{equation}
  \label{105}
  \lim_{t\to\infty}{\V J_t \over \sqrt t} 
 =\frac{1}{4\sqrt\pi}.
\end{equation}
\end{thm}

\proof 
Working from \reff{115}, we get
\begin{eqnarray}
  \nonumber
  \E (J_t)^2 \= \sum_{i< 0} \P(X^{2i}_t\ge 0)\; 
+\; \sum_{i\ge 0} \P(X^{2i}_t< 0) 
 \; +\;2\,\sum_{i<0}\sum_{i<j<0}\, \P(X^{2i}_t\ge 0,\,X^{2j}_t\ge 0) \\
  \label{eq:2} 
&+&\;2\,\sum_{i\ge 0}\sum_{i>j\ge 0}\, \P(X^{2i}_t< 0,\,X^{2j}_t< 0) 
\;-\;2\sum_{i<0}\sum_{j\ge 0}\, \P(X^{2i}_t\ge 0,\,X^{2j}_t< 0)\\ \nonumber\\
  \nonumber 
  (\E J_t)^2 \= \sum_{i< 0} \P^2(X^{2i}_t\ge 0)\; +\; \sum_{i\ge 0}
  \P^2(X^{2i}_t< 0) 
 \; +\;2\,\sum_{i<0}\sum_{i<j<0}\, \P(X^{2i}_t\ge 0)\,\P(X^{2j}_t\ge 0) \\
  \label{eq:3}
&+&\;2\,\sum_{i\ge 0}\sum_{i>j\ge 0}\, \P(X^{2i}_t< 0)\,\P(X^{2j}_t< 0) 
\;-\;2\sum_{i<0}\sum_{j\ge 0}\, \P(X^{2i}_t\ge 0)\,\P(X^{2j}_t< 0)
\end{eqnarray}
Immediately we have:
\begin{equation}
  \label{eq:4}
  \sum_{i< 0} \P(X^{2i}_t\ge 0) + \sum_{i\ge 0} \P(X^{2i}_t< 0) = \sum_{i< 0}
  \P(X^{i}_t\ge 0) 
\end{equation}
and analogously for the $\P^2$ terms in~(\ref{eq:3}).  Using
\begin{equation}
  \label{eq:5}
  \P(A\,B) - \P(A)\,\P(B) \,=\, - (\P(A\,B^c) - \P(A)\,\P(B^c))
  \,=\, \P(A^c\,B^c) - \P(A^c)\,\P(B^c)
\end{equation}
we get
\begin{equation}
  \label{eq:6}
  \V J_t \,=\, {\cal V}_t\, +\, {\cal E}_t,
\end{equation}
where
\begin{equation}
  \label{v0}
  {\cal V}_t=\sum_{i< 0}\P(X^{i}_t\ge 0)-\sum_{i< 0}\P^2(X^{i}_t\ge 0)
\end{equation}
and
\begin{eqnarray}
 \nonumber
 {\cal E}_t  \= \, \Bigl(\sum_{i,j<0,\,i\ne j}\,+\,\sum_{i,j\ge0,\,i\ne j}
                  \,+\,2\sum_{i<0,j\ge0}\Bigr)\Bigl(\P(X^{2i}_t\ge
  0,\,X^{2j}_t\ge 0) - \P(X^{2i}_t\ge  0)\,\P(X^{2j}_t\ge 0)\Bigr)\\
  \label{eq:7}
  \= \, \sum_{i\ne j}\Bigl(\P(X^{2i}_t\ge
  0,\,X^{2j}_t\ge 0) - \P(X^{2i}_t\ge  0)\,\P(X^{2j}_t\ge 0)\Bigr)
\end{eqnarray}

Since $\P(X^{i}_t\ge 0)\P(X^{j}_t\ge 0)= \P(Y^{i}_t\ge 0,Z^{j}_t\ge 0)$, we
can use \reff{110e} with $f(i,j) =\one\{i\ge 0, j\ge 0\}$ to get
\begin{eqnarray}
 \nonumber
\lefteqn{\P(X^{i}_t\ge
  0,\,X^{j}_t\ge 0) - \P(X^{i}_t\ge  0)\,\P(X^{j}_t\ge 0)}\\
 \label{eq:8}
&=&-\frac12 \,\int_0^t \sum_{y} \P(
  \{X^{i}_s, \,X^{j}_s\} = \{y,y +1\} )\, 
   \Bigl(\P(Y^y_{t-s}\ge 0 )-\P(Y^{y+1}_{t-s}\ge 0 )\Bigr)^2 \,ds
\end{eqnarray}
See also Theorem 2 of Ferrari, Galves and Landim (2000) for a probabilistic
proof of the previous identity. Translation invariance and self-duality of
$(X^{i}_s, \,X^{j}_s)$ implies that \reff{eq:8} equals
\begin{eqnarray}
 -\frac12\,\int_0^t \sum_{y} \P(\{X^0_s, \,X^{1}_s\} = \{i-y,j-y\} ) 
   \P^2(Y^0_{t-s}=y)\,ds \label{8}
\end{eqnarray}
{From}~(\ref{eq:7}) and~(\ref{eq:8})--\reff{8}, 
\begin{equation}
 \label{eq:9}
{\cal E}_t\;=\; -\frac12\,\int_0^t \sum_{y}  \P^2(Y^0_{t-s}=y)
                \sum_{i\ne j}\P(\{X^0_s, \,X^{1}_s\} = \{2i-y,2j-y\} )\,ds.
\end{equation}
Since $X^0_t\neq X^1_t$,  
\begin{equation}
 \label{eq:10}
\sum_{i\ne j}\P(\{X^0_s, \,X^{1}_s\} = \{2i-y,2j-y\}) 
\;=\;2\,\P(X^{1}_s(\mod\, 2)=X^0_s(\mod\, 2)= y (\mod\, 2))
\end{equation}
Let $A_y= \{(i,j)\in\Z^2: i(\mod\, 2)=j(\mod\,
2)=y(\mod\, 2)\}$. We show below that 
\begin{equation}
 \label{eq:11}
\lim_{s\to\infty} \P((X^{0}_s,X^1_s)\in A_y) \,=\,
1/4\,.  
\end{equation}
uniformly in $y$. Hence,
\begin{equation}
 \label{eq:15}
 \lim_{t\to\infty}t^{-1/2}{\cal E}_t
 = -\frac14\lim_{t\to\infty}t^{-1/2}\int_0^t\sum_{y}\P^2(Y^0_{s}=y)\,ds
\end{equation}
Let $Z^0_t$ be an independent copy of $Y^0_t$. Since $\sum_{y}\P^2(Y^0_{s}=y)=
\P(Y^0_{t}-Z^0_{t}=0)$, changing variables the above limit equals
\begin{eqnarray}
 \label{16c}
\lim_{t\to\infty}\int_0^1(st)^{1/2}\P(Y^0_{st}-Z^0_{st}=0)\,\frac{ds}{\sqrt s}
&=&\int_0^1\lim_{t\to\infty}(st)^{1/2}\P(Y^0_{st}-Z^0_{st}=0)\,\frac{ds}{\sqrt
  {s}},
\end{eqnarray}
where the interchange of the limit and the integral are guaranteed by the
local central limit theorem for $(Y^0_t-Z^0_t)$, which is a symmetric random
walk of rate $2$. This also implies that \reff{16c} equals
\begin{eqnarray}
 \label{eq:16}
\int_0^1\frac{1}{\sqrt{2\pi}\sqrt2}\,\frac{ds}{\sqrt {s}}\;=\; \frac1{\sqrt \pi}.
\end{eqnarray}
We conclude that
\begin{equation}
 \label{eq:17}
 \lim_{t\to\infty}t^{-1/2}{\cal E}_t
 \;=\;-\frac{1}{4\sqrt\pi}.
\end{equation}

To compute ${\cal V}_t$ notice that
\begin{eqnarray}
  \label{eq:12}
  \sum_{i< 0} \P(X^{i}_t\ge 0)\=\sum_{i>0} \P(X^0_t\ge i)\;=\;\E((X^0_t)^+);\\
  \label{eq:13}
  \sum_{i< 0} \P^2(X^{i}_t\ge 0)\=\sum_{i>0} \P^2(X^0_t\ge i)
  \;=\;\sum_{i>0} \P(Y^0_t\wedge Z^0_t\ge i)\;=\;\E[(Y^0_t\wedge Z^0_t)^+],
\end{eqnarray}
Thus
\begin{equation}
 \label{eq:14}
 \lim_{t\to\infty}t^{-1/2}{\cal V}_t
 \;=\;\E(X^+)-\E[(X\wedge X')^+]\;=\;\frac1{2\sqrt\pi},
\end{equation}
where $X$ and $X'$ are i.i.d.~standard normals.

Finally , substituting~(\ref{eq:17}) and~(\ref{eq:14}) in~(\ref{eq:6}) we get
\reff{105}. \square

\paragraph{Proof of \reff{eq:11}} 
The continuous time Markov chain $(Y^0_t(\mod\, 2), Y^1_t(\mod\, 2))$
converges exponentially fast to the uniform distribution in
$\{(0,1),(1,0),(0,0),(1,1)\}$. This implies that there exist positive
constants $C_1$, $C_2$ such that
\begin{equation}
 \label{11a}
|\P((Y^{i}_s,Z^j_s)\in A_y) \,-\,
1/4| \le C_1 e^{-C_2 t}\,. 
\end{equation}
uniformly in $i,j,y$. 
Writing $f_y(i,j):=\one\{(i,j)\in A_y\}$ and using \reff{110e} we
get 
\begin{eqnarray}
\nonumber
 \lefteqn{ |\P((X^0_t,X^1_t)\in A_y) - \P((Y^0_t,Z^1_t)\in A_y)|}\nonumber\\
&\le&\frac12 \int_0^t \,ds\,\E\Bigl(\one\{|X^0_s-X^1_s|=1\}\\
&&\qquad\times\,
[f_y(Y^{X^{0}_{t-s}}_t,Z^{X^{1}_{t-s}}_t)+
  f_y(Y^{X^{1}_{t-s}}_t,Z^{X^{0}_{t-s}}_t) -
  f_y(Y^{X^{0}_{t-s}}_t,Z^{X^{0}_{t-s}}_t) 
  - f_y(Y^{X^{1}_{t-s}}_t,Z^{X^{1}_{t-s}}_t)]\Bigr)\nonumber\\
\label{11d}
  &\le& 2 \int_0^t \,ds\,\P(|X^1_s-X^0_s|=1)\, C_1\, e^{-C_2 (t-s)}
\end{eqnarray}
(using \reff{11a}).  Now $|X^1_s-X^0_s|$ is a Markov chain in $\{1,2,\ldots\}$
with rates $p(1,2)=p(x,x+1)=p(x,x-1)=1/2,\, x>1$. It can be easily coupled to
a a Markov chain in $\{0,1,2,\ldots\}$ starting in $0$, say $R^0_t$, with
rates $p(0,1)=1,\, p(x,x+1)=p(x,x-1)=1/2,\, x>0$ in such a way that
$|X^1_s-X^0_s|\geq R^0_s$ for all $s$. Since $R^0_t$ is a simple symmetric
random walk reflected at the origin, we get that
$\lim_{s\to\infty}\P(|X^1_s-X^0_s|=1)\leq\lim_{s\to\infty}\P(|R^0_s|\leq1)=0$
and thus, from~\reff{11d} and dominated convergence (after a change of
variables $s\to t-s'$), it follows that $\lim_{t\to\infty}\P((X^0_s,X^1_s) \in
A_y)=\lim_{t\to\infty}\P((Y^0_s,Z^1_s) \in A_y)=1/4$. \square

\begin{lemma}
  \label{1103}
Let $H_t$ be as in \reff{1102} and 
$\tilde H_t=t^{-1/4}[H_t-\E (H_t)].$ Then for all $\la\in\R$
\begin{equation}
\label{eq:expmom}
\limsup_{t\to\infty}\E(e^{\la \tilde H_t})\;\leq \;e^{s\la^2/2},
\end{equation}
where $s=1/\sqrt{2\pi}$.
\end{lemma}

\noindent{\bf Proof.} 
\begin{equation}
\label{eq:l6}
\E(e^{\la \tilde H_t})\;=\;
\frac{\E\left[\exp\left(\la t^{-1/4}\sum_{i<0}\one\{X^{2i}_t\ge0\}\right)
      \right]}
     {\exp\left(\la t^{-1/4}\sum_{i<0}\P(X^{2i}_t\ge0)\right)}
\end{equation}

We will show that the quotient in~(\ref{eq:l6}) is bounded above by a
constant. For that, we need to evaluate the expected value in that equation.
Let $\la\geq0$. 
We will argue below that
\begin{equation}
\label{eq:l7}
\E\left[\exp\left(\la t^{-1/4}\sum_{i<0}\one\{X^{2i}_t\ge0\}\right)\right]
\;\leq\;\prod_{i<0}\E\left[\exp\left(\la t^{-1/4}\one\{X^{2i}_t\ge0\}\right)\right].
\end{equation}

The last expectation equals
\begin{equation}
\label{eq:l8}
1+[\exp(\la t^{-1/4})-1]\P(X^{2i}_t\ge0)= 1+[\la t^{-1/4}+(\la^2 t^{-1/2}/2)
 + o(t^{-1/2})]\P(X^{2i}_t\ge0),
\end{equation}
for all $t$ large enough. The last expression is bounded above by
\begin{equation}
\label{eq:l9}
\exp\left\{[\la t^{-1/4}+(\la^2 t^{-1/2}/2) + o(t^{-1/2})]\P(X^{2i}_t\ge0)\right\}.
\end{equation}
Substituting into the right hand side of~(\ref{eq:l7}), we get 
\begin{equation}
\label{eq:l10}
\exp\left\{[(\la^2 t^{-1/2}/2) + o(t^{-1/2})]\sum_{i<0}\P(X^{2i}_t\ge0)\right\}.
\end{equation}
as an upper bound for the quotient in~(\ref{eq:l6}). It is not difficult 
to see that the expression on the exponent in~(\ref{eq:l10}) converges
to $e^{s\la^2/2}$. 

To finish the argument for $\la\geq0$, we have to justify the 
inequality~(\ref{eq:l7}). That
follows from taking limits as $M\to-\infty$ (and using monotone convergence)
on the respective inequalities 
gotten by replacing the infinite sums by $\sum_{M<i<0}$. These are justified 
by the fact that the functions $\exp(t^{-1/4}\sum_{M<i<0}\one\{X^{2i}_t>0\})$
are bounded, symmetric and positive definite for all $M<0$. The inequalities
then follow
from Proposition 1.7, Chapter VIII of Liggett (1985).

For the case $\la<0$, we use the identity
\begin{equation}
\label{eq:l1}
\sum_{i< 0} \one\{X^{2i}_t\ge 0\} -\sum_{i< 0} \P(X^{2i}_t\ge 0)
=-\left[\sum_{i< 0} \one\{X^{2i}_t< 0\} -\sum_{i< 0} \P(X^{2i}_t< 0)\right]
\end{equation}
and a similar argument as above. \square

\begin{lemma}
  \label{110}
  Let $\eta_0$ be given by the flat condition as in Theorem \ref{104}. Then
\begin{equation}
\label{eq:expabs}
\sup_{t\geq0}\E(e^{|J_t/t^{1/4}|})<\infty.
\end{equation}
Furthermore for all $K>0$ and all $t$ large enough
\begin{equation}
\label{eq:lema}
\P(|J_t|>Kt^{1/4}\log t)\;\leq \;c t^{-K},
\end{equation}
where $c$ is a constant. 
\end{lemma}

\noindent{\bf Proof.} The bound \reff{eq:expabs} follows straightforwardly
from Lemma \reff{1103} and relations \reff{1102}. 

{From} the relations \reff{1102}, to show \reff{eq:lema} it is enough to prove
the result with $|H_t-\E (H_t)|$ replacing $|J_t|$ in \reff{eq:lema} (the
constant of course does not need to be the same).  We have
\begin{equation}
\label{1105}
\P(|H_t|>Kt^{1/4}\log t)=\P(|\tilde H_t|>\log t^K)\leq c' t^{-K}
\end{equation}
where the last inequality follows from the exponential Markov inequality and
$c'=\sup_{t\geq0}\E(e^{|\tilde H_t|})$ is finite by \reff{eq:expabs}.  \square

\noindent{\bf Graphical construction of free process.}
 Let $\zeta_t$ be the
process defined by
\begin{equation}
  \label{112}
  d \zeta_t(x) = \Delta\zeta_t(x)\, dN_t(x),  
\end{equation}
where the discrete Laplacian $\Delta$ was defined in \reff{p28}. In words,
each time a Poisson mark of the process $N_t(x)$ occurs, the height at $x$ at
time $t$ decreases or increases two unities, according to the value of the
Laplacian at this point at this time; if the Laplacian vanishes, no jump
occurs. This process has generator \reff{gen0}.

\begin{lemma}
  \label{103}
Let $\eta_0(x)  = \zeta_0(x+1) -\zeta_0(x)$. Then 
\begin{equation}
  \label{113}
  \eta_t(x) = \zeta_t(x+1)-\zeta_t(x) 
\end{equation}
where the
processes $\zeta_t$ and $\eta_t$ are defined by \reff{111} and \reff{112} and
have initial conditions $\eta_0$ and
$\zeta_0$ respectively. Furthermore,
\begin{equation}
  \label{102}
  \zeta_t (0) -\zeta_0(0) = 2J_t\,.
\end{equation}
\end{lemma}

\proof Notice that from \reff{113}, 
\begin{eqnarray}
  \label{x1}
\Delta\zeta_t(0) =2 (\eta_t(-1)- \eta_t(0))
\end{eqnarray}
Assume that there is a mark of the process $N_t(-1)$ at time $t$. Then
\reff{x1}, \reff{111} and \reff{112} imply that if $\eta_t(-1)- \eta_t(0)=0$
no changes occur neither for $\eta_t(-1),\eta_t(0)$ nor for $\zeta_t(0)$; if
$\eta_t(-1)- \eta_t(0)=1$, an exclusion particle jumps from $-1$ to $0$ and
the free process at the origin jumps two units up; if $\eta_t(-1)-
\eta_t(0)=-1$, an exclusion particle jumps from $0$ to $-1$ and the free
process at the origin jumps two units down. Identity \reff{102} follows from
\reff{113}. \square

\begin{lemma}
  \label{110q}
Let $\zeta_t$ be the free process with flat initial condition \reff{99}. Then  
\begin{equation}
  \label{105p}
  \lim_{t\to\infty}{\V \ze_t(0) \over \sqrt t} 
 =\frac{1}{\sqrt\pi};
\end{equation}
\begin{equation}
\label{110r}
\sup_{t\geq0}\E(e^{|\ze_t(0)/t^{1/4}|})<\infty
\end{equation}
and for all $K>0$ and all $t$ large enough
\begin{equation}
\label{110p}
\P(|\ze_t(0)|>Kt^{1/4}\log t)\;\leq \;c t^{-K}\,.
\end{equation}
\end{lemma}
\proof It follows from identity \reff{102}, the limit \reff{105} and the
bounds \reff {eq:expabs} and \reff{eq:lema}. \square

\section{Coupling the wall and the free processes}

We construct graphically the wall process which simultaneously provides
another graphical construction for the free process. Under this construction
the wall process dominates the free one. We consider two independent families
of Poisson processes with the same law as $N_t(x)$ called $N^+_t(x)$ and
$N^-_t(x)$, to be used for upwards and downwards jumps, respectively. The
process satisfying the equations
\begin{equation}
  \label{121}
  d \xi_t(x) \;=\; \Delta
  \xi_t(x)\,\one\{\Delta
  \xi_t(x)  > 0\}\, dN^+_t(x)\;+\;\Delta
  \xi_t(x)\,\one\{\Delta
  \xi_t(x)  < 0,\,\xi_t(x) + \Delta
  \xi_t(x)  \ge 0\}\, dN^-_t(x)
\end{equation}
has generator \reff{gen}. The process $\zeta_t$ satisfying
\begin{equation}
  \label{121a}
  d \ze_t(x) = \Delta
  \xi_t(x)\,\one\{\Delta
  \xi_t(x)  > 0\}\, dN^+_t(x)+\Delta
  \xi_t(x)\,\one\{\Delta
  \xi_t(x)  < 0\}\, dN^-_t(x).  
\end{equation}
has generator \reff{gen0}.

In words, when a time event of the process $N^+_t(x)$ occurs at time $t$, the
process $\xi_t$ at site $x$ and time $t$ jumps two units upwards if $\Delta
\xi_t(x)>0$. When a time event of the process $N^-_t(x)$ occurs at time $t$,
the process $\xi_t$ at site $x$ and time $t$ jumps two units downwards if
$\Delta \xi_t(x)<0$ \emph{and} the wall condition $\xi_t(x) + \Delta \xi_t(x)
\ge 0$ holds. The process satisfying \reff{121a} follows the same marks in the
same manner but ignoring the wall condition. The difference with the process
satisfying \reff{112} is that in this case the Poisson events $N_t$ are used
for \emph{both} upwards and downwards jumps; this construction is not
attractive in the sense that it does not satisfy \reff{p12} below.

Let $r$ be a non negative integer and
$\xi^r_t$ and $\zeta^r_t$ be the processes defined by \reff{121} and
\reff{121a} but with initial condition
\begin{equation}
  \label{p13}
  \xi^r_0(x) \,=\, \zeta^r_0(x) \,=\, r + x(\mod 2)
\end{equation}
Notice that $\ze^0_t$ and $\ze_t$ as defined in \reff{112} have the same law
but are \emph{different} processes. The processes $\xi_t$ and $\xi^r_t$
defined by \reff{121} and the same initial condition satisfy
\begin{equation}
  \label{p12}
  \xi_t(0) \le \xi^r_t(0)
\end{equation}
for all $r\ge 0$. This joint construction corresponds to what Liggett (1985)
calls \emph{basic coupling}.

\begin{lemma}
  \label{p10}
There exists a constant $c>0$ such that for any $K>0$ and $t\geq0$
  \begin{equation}
    \label{p11}
    \P(\xi_t(0) > 2K t^{1/4}\log t) \le c t^{2-K}
  \end{equation}
\end{lemma}

\proof Let $a_t = 2K t^{1/4}\log t$. Take $r\ge 0$ and write
\begin{eqnarray}
  \label{p14}
  \P(\xi_t(0)\ge a_t) &\le& \P(\xi^r_t(0)\ge a_t) \nonumber\\
&=& \P(\xi^r_t(0)\ge a_t\,,\,\xi^r_t(0) = \zeta^r_t(0))\,+\, \P(\xi^r_t(0)\ge
a_t\,,\,\xi^r_t(0) \neq \zeta^r_t(0) )\nonumber\\ 
&\le& \P(\zeta^r_t(0)\ge a_t)\,+\, \P(\xi^r_t(0) \neq \zeta^r_t(0) )
\end{eqnarray}
The first term in \reff{p14} will be bounded using Corollary \ref{110q}. To
bound the second term notice that if the interacting process and the free
process differ at the origin this is due to a collision of the interacting
process with the wall at some point $x$ that separate the two processes at $x$
at some time $s$; the discrepancy then propagates and arrives to zero by time
$t$. We fix an $\alpha>0$ and separate the discrepancies in two classes: those
that come from the interval $[-\alpha t,\alpha t]$ and those that come from
outside this interval. If in the time interval $[0,t]$ the free process does
not touch the wall in the space interval $[-\alpha t,\alpha t]$, then the
discrepancy must come from outside. Hence,
\begin{eqnarray}
\nonumber
  \{\xi^r_t(0)\neq \zeta^r_t(0) \}\!\!\!\!&\subset&\!\!\!\! 
\{\zeta^r_t(x)<0\,\hbox{ for some }s\in[0,t],\, x\in[-\alpha t,\alpha t]\}\\
\label{p15}
&\,\,\,\,\,\,\,\,\,\,\,\,\cup&\!\!\!\! \{\hbox{ a discrepancy from }[-\alpha
  t,\alpha t]^c\hbox{ reaches } 0 \hbox{ up to time }t\}
\end{eqnarray}
Observe that the law of $\zeta^r_t(x)-r$ is the same as the law of $\ze_t(0)$
and that $\P(\ze_s(0) <-r)\le \P(\ze_s(0) > r)$ due to 
the initial condition being non-negative. Hence, fixing
\begin{equation}
  \label{1234}
  r= a_t/2\,,
\end{equation} 
the probability of the first
event in the right hand side of \reff{p15} is bounded by
\begin{eqnarray}
  \label{p23}
   (2\alpha t+1)\, \P(\ze_s(0) > a_t/2\, \hbox{ for
  some } s\in[0,t])
\end{eqnarray}
{From} \reff{110r} and the exponential Markov inequality, we have that
$\sup_{s\leq t}\P(\ze_s(0)>r)\le ct^{-K}$ for some constant $c$, so we can
bound \reff{p23} with
\begin{equation}
  \label{p24}
  (2\alpha t+1)\, \int_0^t c t^{-K} ds \;\le\; (2\alpha+1) \, c\, t^{2-K}
\end{equation}

To bound the probability of the second event in the right hand side of
\reff{p15} notice that discrepancies cannot travel faster than $N_t$, a Poisson
process of parameter 1. Hence
\begin{eqnarray}
  \label{p17}
  \P(\hbox{a discrepancy from }[-\alpha
  t,\alpha t]^c\hbox{ reaches } 0 \hbox{ up to time }t) \le
  2\P(N_t > \alpha t)\le 2\,e^{-t (\alpha+1-e)}
\end{eqnarray}
using the exponential Chebyshev inequality. Fixing $\alpha=2$, and using the
bounds \reff{p23} and \reff{p17}, the probability of \reff{p15} is bounded by
\begin{equation}
  \label{p31}
 4 \, c\, t^{2-K} \,+\,2\,e^{-t (3-e)}\le
  c't^{2-K} 
\end{equation} 
for some constant $c'$ and sufficiently large $t$. \square

\medskip
\noindent{\bf Proof of Theorem \ref{98}.}
It follows straightforwardly {from} (\ref{110r}) that 
${\tilde\ze_t^2}$ is
uniformly integrable, where ${\tilde\ze_t}=t^{-1/4}\ze_t(0)$.
This, together with \reff{105p}, implies the lower
bound in \reff{100}, as we will see now. Indeed,
\reff{cou} implies that $t^{-1/4}\E\xi_t(0)\geq\E|\tilde\ze_t|$.
Now,
\begin{eqnarray}
  \label{p25}
  \V\tilde\ze_t\leq\E\tilde\ze_t^2=
\E(\tilde\ze_t^2;\,\tilde\ze_t^2\leq M^2)
+\E(\tilde\ze_t^2;\,\tilde\ze_t^2>M^2)
\leq M\E|\tilde\ze_t|+\e_M
\end{eqnarray}
uniformly in $t$, where $M$ is an arbitrary positive number and
$\e_M\to0$ as $M\to\infty$. Thus 
$t^{-1/4}\E\xi_t(0)\geq(\V\tilde\ze_t-\e_M)/M$. We conclude that
\begin{eqnarray}
  \label{p26}
\liminf_{t\to\infty}t^{-1/4}\E\xi_t(0)
\geq\sup_{M>0}\left(\frac1{\sqrt\pi}-\e_M\right)/M>0.
\end{eqnarray}

For the upperbound, we use \ref{110p} to obtain
\begin{eqnarray}
  \label{p27}
  {\E \xi_t(0) \over t^{1/4}\log t} = \sum_{k\ge 0} \P(\xi_t(0) > k
  t^{1/4}\log t)
\le 4+\sum_{k\ge 5}c't^{2-k/2} \;\le c_2 \;<\;\infty
\end{eqnarray}
for some constant $c_2<\infty$. \square

\section{Numerical simulation}

We have simulated the processes $\zeta_t$ and $\xi_t$ numerically, using
various pseudo-random number generators.
The interface is of length $L$ with periodic boundary conditions,
so that the processes live in
\begin{equation}
\X_{0L}=\{\zeta\in\Z^{\Z/L\Z}: |\zeta(x)- \zeta(x+1)|=1,\, \zeta(0) 
\hbox{ even}\} 
\end{equation}
or
\begin{equation}
\X_L=\{\xi\in\N^{\Z/L\Z}: |\xi(x)- \xi(x+1)|=1,\, \xi(0) \hbox{ even}\} 
\end{equation}
Time is an integer multiple of $2L^{-1}$, i.e. $t\in (2L^{-1})\N$.
For each time step, a site $x$ is chosen randomly according to the uniform 
measure on $\Z/L\Z$, and the interface is updated with the same rules
as in the continuous time version of the processes. 
The transition operator corresponding to a time step $\delta t=2L^{-1}$ is
\begin{equation}
T_0f(\zeta)=f(\zeta)+L^{-1}\sum_x \, [f(\zeta + \Delta
  \zeta(x)\,\delta_x) - f(\zeta)]
\end{equation}
or
\begin{equation}
Tf(\xi)=f(\xi)+ L^{-1}\sum_x \one\{\xi + \Delta
  \xi(x) \,\delta_x \ge 0\}\, [f(\xi + \Delta
  \xi(x)\,\delta_x) - f(\xi)]
\end{equation}
  $Tf(\xi)$ is the expected value of the function $f$ evaluated at time
  $2L^{-1}$ (after one step) when the initial configuration is $\xi$ for the
  discretized version of the process. The same interpretation is valid for
  $T_0f(\zeta)$. By abuse of notation we call the
  discrete time versions of the process $\xi_t$ and $\zeta_t$ as we did for
  the continuous time versions. Notice that 
\begin{equation}
\lim_{L\to\infty} {Tf(\xi)-f(\xi))\over\delta t} = {\cal L}f(\xi)\,;\qquad
\lim_{L\to\infty} {T_0f(\zeta)-f(\zeta))\over\delta t} = {\cal L}_0f(\zeta)\,.
\end{equation}
The processes $\xi_t$ and $\zeta_t$ are coupled in the simplest possible
way: the same random sequence of sites are used for both.  Notice however that
this coupling is different from the one described in Section 3 (in particular
it is not attractive in the sense that it does not necessarily satisfy
\reff{cou} but it is faster). For $L$ finite the discrete time and
continuous time processes can be identified up to a time change, using the
ordered sequence of updated sites. The random time involved in the time change
has fluctuations which should be negligible for our purposes.  

The numerical samples for the data shown below were drawn using either the
Mersenne Twister pseudorandom integer generator, see
Matsumoto,  Nishimura (1998), or the R250 pseudorandom generator, see
Kirkpatrick, Stoll (1981). The length $L$ is $10^6$ or $2^{20}=1024^2$ and
time runs up to $2.10^6$ or $2^{21}$.
The number of calls to the generator for the realization of one sample
of length $L$ up to time $t$ is $L.t/2\le10^{12}$, which of course is much
less than the period of the generator (a necessary but not sufficient
condition for reliability). 
We compute empirical averages
\begin{equation}
\overline{\xi_t^2}=L^{-1}\sum_x\xi_t(x)^2\ ,\quad 
\overline{\zeta_t^2}=L^{-1}\sum_x\zeta_t(x)^2
\end{equation}
and empirical distribution functions 
\begin{equation}
f_t(n)=L^{-1}\sum_x\one\{\xi_t(x)=n\}\ ,\quad 
f_{0,t}(n)=L^{-1}\sum_x\one\{\zeta_t(x)=n\}\ ,\quad n\in\Z
\end{equation}
scaled into
\begin{equation}
\phi_t(s)=t^{1/4}L^{-1}\sum_x\one\{t^{-1/4}\xi_t(x)=s\}\ ,\quad 
\phi_{0,t}(s)=t^{1/4}L^{-1}\sum_x\one\{t^{-1/4}\zeta_t(x)=s\}\ ,\quad 
s\in t^{-1/4}\Z
\end{equation}
which, extended to $s\in\R$ approximate the Schwartz distributions
\begin{equation}
\tilde\phi_t(s)=L^{-1}\sum_x\delta(t^{-1/4}\xi_t(x)-s)\ ,\quad 
\tilde\phi_{0,t}(s)=L^{-1}\sum_x\delta(t^{-1/4}\zeta_t(x)-s)\ ,\quad 
s\in\R
\end{equation}
where $L^{-1}\sum_x$ is an ersatz for the expectation over a real random
variable, limit of $t^{-1/4}\xi_t(x)$.

The processes were studied for time $t\le 2L$,  whereas the effect of finite
size with periodic boundary conditions is expected to be visible only after a 
time of order $L^2$, the relaxation time of an interface of length $L$.
The law of large numbers in empirical averages as above is believed to be at
work with an effective number of weakly dependent variables of order
$L/t^{1/2}$: the interface at time $t$ can be thought of as a collection of 
$L/t^{1/2}$ segments of
length $t^{1/2}$, the different segments being weakly dependent.
For $t=L$ and one sample of the process, we have only $t^{1/2}$ 
independent segments, hence an expected relative statistical error of order 
$t^{-1/4}$. This explains the more erratic behaviour at larger times
in Fig. 1.

\vspace{0.3cm}
 {\centering \resizebox*{0.9\textwidth} {!}{\includegraphics{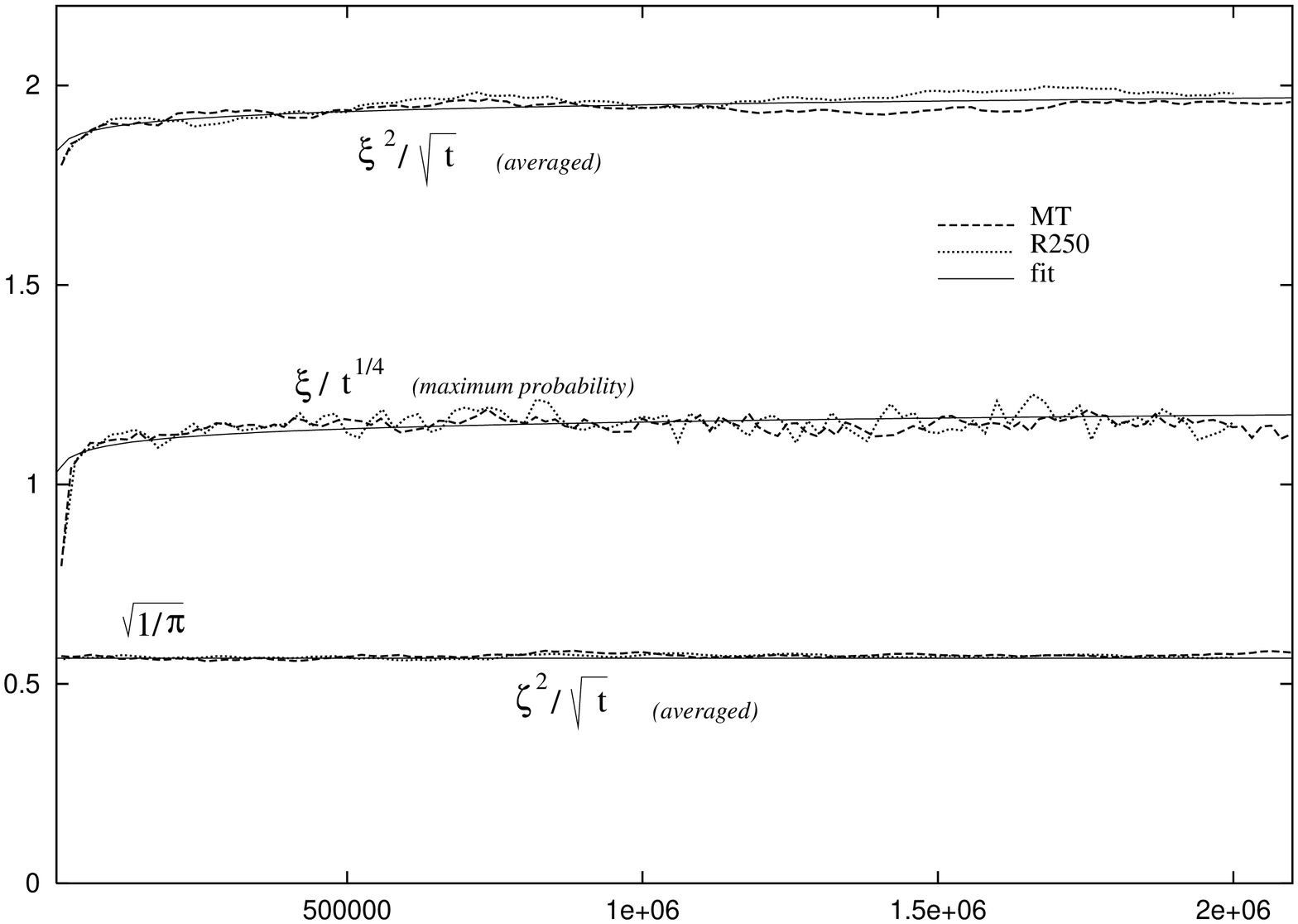}} 
 \par \footnotesize
 Fig. 1: from top to bottom, as function of time : 
 $t^{-1/2}\overline{\xi_t^2}$ together with best fit
 $1.62+0.024\log t$ ;
 the value of $s$ where $\phi_t(s)$ is maximum, together with best fit
 $\sqrt{0.55+0.057\log t}$~; and
 $t^{-1/2}\overline{\zeta_t^2}$ together with the exact asymptotic value 
 $\sqrt{1/\pi}$.
 Graphs labelled MT are averages over 6 runs with the MT random generator with 
 different seeds. Graphs labelled R250 are averages over 5 runs with the R250 
 random generator with different seeds. Interface length is $L=2^{20}$ or 
 $L=10^6$.}
\vspace{0.3cm}

The numerical experiment clearly favors an asymtotic behavior
$\E\xi_t(0)^2\sim t^{1/2}\log t$ as $t\to\infty$.

\vspace{0.3cm}
{\centering \resizebox*{0.9\textwidth} {!}{\includegraphics{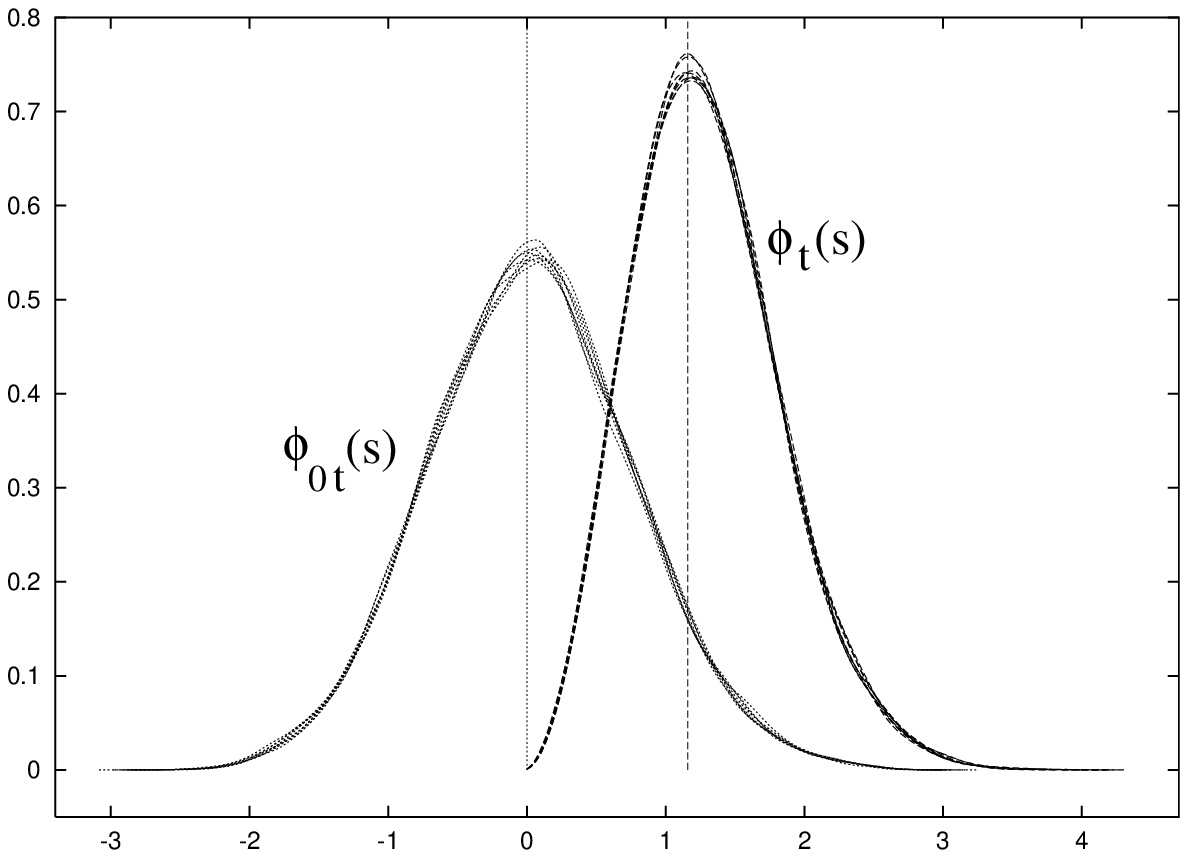}} 
 \par  
 \footnotesize 
 \hfil Fig.~2: $\phi_t(s)$ and  $\phi_{0,t}(s)$ at large times, from the same 
 data as in Fig. 1.
 }
\vspace{0.3cm}

Figure 2 shows the scaled empirical distribution functions at various large 
times. Clearly $t^{-1/4}\zeta_t(x)$ converges to a centered Gaussian
random variable as expected. The distributon function of $t^{-1/4}\xi_t(x)$ 
is markedly asymmetrical. Zooming around $s=0$ indicates
$\phi_\infty(0)=\phi_\infty'(0)=0$ and $\phi_\infty''(0)>0$,
and 
\begin{equation}
f_t(0)=L^{-1}\sum_x\one\{\xi_t(x)=0\}\sim t^{-1/2}
\end{equation}

\section{Interface of the Ising model at zero temperature}\label{ising}
In this section we explain the relation of our model with the interface of a
particular Ising model at zero temperature. Let the ``inverse temperature''
$\beta\ge 0$ and $\sigma_t\in\{-1,+1\}^{\Z^2}$ be the Ising model with
generator
\[
{\cal L}_{\beta}f(\sigma) = \frac12\sum_{x\in\Z^2}
c_\beta(x,\sigma)\, [f(\sigma^x)-f(\sigma)]
\]
with $\sigma^x (z) = \sigma (z)$ for $z\neq x\in\Z^2$, $\sigma^x(x) =
-\sigma(x)$ and $c_\beta(x,\sigma)$ are the Glauber rates
\[
c_\beta(x,\sigma) = \frac{e^{-\beta H(\sigma^x)}}{ e^{-\beta H(\sigma^x)}+e^{-\beta H(\sigma)}}
\]
with Hamiltonian 
\[
H(\sigma) =- \sum_x\sum_{y:|y-x|=1}\sigma(x)\sigma(y) - h \sum_{x:x_1> x_2}
\sigma(x) 
\]
for some positive magnetic field $h>0$. Consider the case $\beta=\infty$ and
assume that the starting configuration $\sigma_0$ is ``all ones'' below
the diagonal and ``all minus ones'' above or in the diagonal:
\begin{equation}
  \label{s00}
\sigma_0(x_1,x_2) = \cases{ +1, &if $x_1>x_2$; \cr
-1, &if $x_1\leq x_2$.}
\end{equation}
In this case for all $t$ the configuration $\sigma_t$ has the property that
all sites have either exactly two or no neighbor with opposite sign;
furthermore, only sites above or in the diagonal may be negative. As a
consequence, for all $t\ge 0$ the rates $c(x,\sigma_t)$ are positive only for
sites $x$ above or in the diagonal for which there are exacly two neighboring
sites with different sign: under initial condition \reff{s00},
\[
c_\infty(x,\sigma_t) = \cases {
1/2 &if $\sum_{y:|y-x|=1}
  \one\{\sigma_t(y)\neq\sigma_t(x)\}=2$ and $x_1\le x_2$\cr
0 &otherwise}
\]
To get the wall process of Theorem \ref{98} from the above dynamics with
initial condition \reff{s00}, we first rotate the lattice by $-45^0$ and
multiply by $\sqrt 2$, that is, we perform the transformation
$R:\Z^2\to{\Z}^2_2;R(x,y)=(x+y,x-y)$, where
${\Z}^2_2:=\{(x,y)\in\Z^2:x+y\mbox{ is even}\}$ is the even sublattice of
$\Z^2$. The above dynamics then induces a dynamics in $\{-1,+1\}^{\Z^2_2}$
given by $\tilde\sigma_t(z)=\sigma_t(R^{-1}z)$, $z\in{\Z}^2_2,t\geq0$.
Defining $\tilde\xi_t(x):= \min\{y:(x,y)\in{\Z}^2_2\mbox{ and
  }\tilde\sigma_t(x,y)=-1\}$, $x\in\Z,t\geq0$, we have that the wall process
$\xi_\cdot(\cdot)$ with generator \reff{gen} and initial configuration
\reff{99} has the same law as $\tilde\xi_\cdot(\cdot)$.

{\thebibliography{}
\parskip 0mm

\item D.B. Abraham, P. Upton (1989) 
Dynamics of Gaussian interface models.
{\sl Phys. Rev B \bf 39}, 736.

\item D.B. Abraham, P. Collet, J. De Coninck, F. Dunlop (1990)
Langevin Dynamics of an Interface near a Wall.
{\sl J. Stat. Phys. \bf 42}, 509-532.

\item E.V. Albano, K. Binder, D.W. Heermann, W. Paul (1992)
{\sl Physica A \bf 183}, 130

\item R. Arratia (1983) The motion of a tagged particle in the simple 
symmetric exclusion system on ${\Z}$.
{\sl Ann. Probab. \bf 11}, no. 2, 362--373. 

\item N. C. Bartelt, J. L. Goldberg, T. L. Einstein, Ellen D. Williams, 
J. C. Heyraud and J. J. M\'etois (1993)
Brownian motion of steps on Si(111).
{\sl Phys. Rev B \bf 48}, 15453-15456.

\item K. Binder (1990) Growth kinetics of wetting layers at surfaces.
pp 31-44 in {\sl Kinetics of Ordering and Growth at Surfaces}. Edited by 
M.G. Lagally. Plenum Press, New-York.

\item B. Blagojevic, P.M. Duxbury (1999) 
Atomic diffusion, step relaxation, and step fluctuations.
{\sl Phys. Rev E \bf 60}, 1279-1291.

\item E. Bolthausen, J.-D. Deuschel, O. Zeitouni (1995) Entropic repulsion of
  the lattice free field. {\sl Comm. Math. Phys. \bf 170}, no. 2, 417--443.

\item J. Bricmont (1990) Random surfaces in statistical mechanics (1990) in
  Wetting phenomena.  {\sl Proceedings of the Second Workshop held at the
    University of Mons}, Mons, October 17--19, 1988. Edited by J.  De Coninck
  and F. Dunlop.  {\sl Lecture Notes in Physics, \bf 354}.  Springer-Verlag,
  Berlin.
 
\item J. Bricmont, A. El Mellouki, J. Fr\"ohlich (1986) Random surfaces in
  Statistical Mechanics - roughening, rounding, wetting {\sl J. Stat.
    Phys. \bf 42}: (5-6) 743-798

\item  F. Cesi,  F. Martinelli (1996) On the layering
  transition of an SOS surface interacting with a wall. I. Equilibrium
  results. {\sl J. Stat.  Phys. \bf 82}, no. 3-4, 823--913.
  
\item J. De Coninck, F. Dunlop, F. Menu (1993) Spreading of a Solid-On-Solid
  drop.  {\sl Phys. Rev. E \bf 47}: (3) 1820-1823.
 
\item A. De Masi, P. A. Ferrari (2002) Flux fluctuations in the one dimensional
  nearest neighbors symmetric simple exclusion process. To appear in {\sl
    J. Stat.  Phys.} http://front.math.ucdavis.edu/math.PR/0103233. 

\item J.-D. Deuschel (1996) Entropic repulsion of the lattice
free field. II. The $0$-boundary case. {\sl Comm. Math. Phys. \bf 181}, no. 3,
  647--665.
  
\item J.-D. Deuschel, G. Giacomin (1999) Entropic repulsion for the free
  field: pathwise characterization in $d\geq3$. {\sl Comm.  Math. Phys. \bf
    206}, no. 2, 447--462.
  
\item E. Dinaburg,  A.E. Mazel (1994) Layering transition in SOS model with
  external magnetic-field.  {\sl J. Stat.  Phys. \bf 74}: (3-4)
  533-563.
 
\item P. A. Ferrari, S. Mart\'{\i}nez (1998) Hamiltonians on random walk
  trajectories. {\sl Stochastic Process. Appl. \bf 78}, no. 1, 47--68.

\item P. A. Ferrari, A. Galves, C. Landim (2000)
Rate of convergence to equilibrium of symmetric simple exclusion processes. 
{\sl Markov Processes and Related Fields \bf 6}, 73-88. 

\item T. Funaki, S. Olla (2001) Fluctuations for $\nabla\phi$ interface model
  on a wall. Stochastic Process. Appl. 94, no. 1, 
1--27.

\item P. Holick\'y, M. Zahradn\'{\i}k (1993) On entropic repulsion in low
  temperature Ising models.  {\sl Cellular automata and cooperative systems}
  (Les Houches, 1992), 275--287, NATO Adv. Sci. Inst. Ser. C Math.  Phys.
  Sci., 396, Kluwer Acad. Publ., Dordrecht.

\item   W. D. Kaigh (1976)
An invariance principle for random walk conditioned by a late return to zero. 
{\sl Ann. Probability \bf 4}, no. 1, 115--121. 

\item  S. Kirkpatrick, E. Stoll (1981)  {\sl J. Computational Physics}
{\bf 40}, 517-526.  

\item J.L. Lebowitz, A.E. Mazel, Y.M.  Suhov (1996) An Ising interface between
  two walls: Competition between two tendencies. {\sl Reviews In Mathematical
    Physics \bf 8}: (5) 669-687.

\item J.L. Lebowitz, C. Maes
(1987) The effect of an external field on an 
  interface, entropic repulsion. {\sl J. Stat. Phys. \bf 46}, no.  1-2,
  39--49.

\item T.M. Liggett (1985) {\sl Interacting particle systems. 
Grundlehren der Mathematischen Wissenschaften
[Fundamental Principles of Mathematical Science], 276.} 
Springer-Verlag, New York-Berlin.
 
\item R. Lipowsky (1985) Nonlinear growth of wetting layers.
{\sl J. Phys. A \bf 18}, L585-L590

\item M. Matsumoto, T. Nishimura (1998) 
Mersenne Twister: A 623-dimensionally 
equidistributed uniform pseudorandom number generator, 
ACM Trans. on Modeling and Computer Simulation {\bf 8}, 3-30. 
(http://www.math.keio.ac.jp/~matumoto/emt.html).
   
\item K. K. Mon, K. Binder, and D. P. Landau (1987)
Monte Carlo simulation of the growth of wetting layers
{\sl Phys. Rev.} B {\bf 35}: (7) 3683 -- 3685.

\item D. Stauffer and D. P. Landau (1989)
Interface growth in a two-dimensional Ising model.
{\sl Phys. rev. B \bf39}, 9650-9651.

}

\vskip 10mm
\baselineskip 14pt
\obeylines
Fran\c cois Dunlop
Laboratoire de Physique Th\'eorique et Mod\'elisation (CNRS-ESA8089)
Universit\'e de Cergy-Pontoise, 95031 Cergy-Pontoise, France 
{\tt dunlop@ptm.u-cergy.fr}
http://www.u-cergy.fr/rech/labo/equipes/ptm/Dunlop
\vskip 5mm
Pablo A. Ferrari, Luiz Renato G. Fontes
IME--USP, Caixa Postal 66281, 
05315-970 - S\~{a}o Paulo, Brazil
{\tt pablo@ime.usp.br, lrenato@ime.usp.br}
http://www.ime.usp.br/$\widetilde{\phantom m} $pablo,
http://www.ime.usp.br/$\widetilde{\phantom m} $lrenato

\end{document}